\documentclass[preprint]{elsarticle}
\usepackage{amsfonts,anysize,amsmath,indentfirst,geometry,xcolor}
\marginsize{30mm}{30mm}{16mm}{16mm}
\usepackage[displaymath,pagewise]{lineno}
\usepackage[dvipdfm,colorlinks,linkcolor=blue,citecolor=blue]{hyperref}
\usepackage{mathrsfs}
\usepackage{amssymb}
\usepackage{array}
\usepackage{graphicx}
\usepackage{booktabs}
\usepackage{epstopdf}
\usepackage{multirow}
\usepackage{algorithm}
\usepackage{algorithmic}
\numberwithin{equation}{section}% number the equations according to Section.
\newtheorem{thm}{Theorem}[section]%number the theorem according to Section.
\newtheorem{lem}[thm]{Lemma}%number the theorem and lemma together.
\newtheorem{definition}{Definition}[section]%number the definition according to Section.
\newdefinition{rmk}{Remark}[section]%number the theorem according to Section.
\newdefinition{cor}{Corollary}[section]%number the theorem according to Section.
\newdefinition{example}{Example}[section]%number the theorem according to Section.
\newproof{pf}{\textbf{Proof}}
\newproof{pot}{Proof of Theorem \ref{thm2}}
\modulolinenumbers[1]
\journal{arXiv}

\begin{document}

\begin{frontmatter}

\title{Perturbation analysis for the periodic generalized
coupled Sylvester equation}

\author[1]{Hanyu Li \corref{cor1}}
\ead{hyli@cqu.edu.cn, lihy.hy@gmail.com}
\author[1]{Shaoxin Wang}
\ead{shaoxin.w@gmail.com}
\author[1]{Chan Zheng}
\ead{zchan1988£À163.com}

\cortext[cor1]{Corresponding author}
\address[1]{College of Mathematics and Statistics, Chongqing University,  Chongqing, 401331, P. R. China }
%\address[2]{Department of Mathematics and Statistics, Auburn University, Alabama, 36849, USA }

\begin{abstract}
In this paper, we consider the perturbation analysis for the periodic generalized coupled Sylvester (PGCS) equation. The normwise backward error for this equation is first obtained. Then, we present its normwise and componentwise perturbation bounds, from which the normwise and effective condition numbers are derived. Moreover, the mixed and componentwise condition numbers for the PGCS equation are also given. To estimate these condition numbers with high reliability, the probabilistic spectral norm estimator and the statistical condition estimation method are applied. The obtained results are illustrated by numerical examples.
\end{abstract}


\end{frontmatter}

%\linenumbers

\section{Introduction}
In this paper, we consider the following matrix equation:
\begin{equation}\label{1.1}
\left\{
  \begin{array}{ll}
    A_kX_k-Y_kB_k = E_k, & \hbox{} \\
     & \hbox{$k=1,\cdots,p,$} \\
    C_kX_{k+1}-Y_kD_k = F_k, & \hbox{}
  \end{array}
\right.
\end{equation}
where $A_{k}$, $C_k\in \mathbb{R}^{m\times m}$, $B_{k}$, $D_k\in \mathbb{R}^{n\times n}$, and $E_{k}$, $F_k\in \mathbb{R}^{m\times n}$ are the given coefficient matrices, and $X_{k}$, $Y_k\in \mathbb{R}^{m\times n}$ are the unknown matrices satisfying $X_{p+1}=X_{1}$. Hereafter,  $\mathbb{R}^{m \times n} $ denotes the set of $m \times n$ real matrices.

The equation \eqref{1.1} is called the periodic generalized coupled Sylvester (PGCS) equation with period $p$ (see e.g.,\cite{Chen12,Hajarian14}). It is easy to find that if $p=1$, the PGCS equation reduces to the generalized coupled Sylvester (GCS) equation, which plays an important role in the linear control systems (see e.g., \cite{Datta04,Kons03}). One of the significant applications of this equation originates from computing the stable eigendecompositions of matrix pencils \cite{Demmel87}. Some numerical methods were provided to compute the solution of the GCS equation (see e.g., \cite{Ding05,Jons02a,Jons02b}). Considering the specific structure of this equation, K{\aa}gstr\"{o}m \cite{Kags94} investigated its perturbation analysis, and derived the normwise backward error, normwise perturbation bounds, and normwise condition number. The derived results generalized the corresponding ones for the classic Sylvester equation given in \cite{Higham93}. Since the normwise condition number cannot accurately reflect the influence of perturbations for some small entries in the data and ignores the structures of both input and output data with respect to scaling, Lin and Wei \cite{Linw07} presented the mixed and componentwise condition numbers for the GCS equation. These two condition numbers were named by Gohberg and Koltracht \cite{Gohb93}. The former measures the errors in output using norms and the input perturbations componentwise, and the latter measures both the errors in output and the perturbations in input componentwise. To estimate the normwise, mixed and componentwise condition numbers for the GCS equation effectively, Diao et al. \cite{Diao13} applied the statistical condition estimation (SCE) method, which was proposed by Kenney and Laub in \cite{Kenney94} and found applications in estimating the condition numbers of linear systems, least squares problem, eigenvalue problem, and matrix equations (see e.g., \cite{Diao13,Diao12,Gudmun97,Kenney98a,Keney98b,laub08}). Moreover, the authors also derived the effective condition numbers for the GCS equation and the classic Sylvester equation in \cite{Diao13}, which can be much tighter than the normwise ones in \cite{Higham93,Kags94} in some cases.

The PGCS equation also finds applications in many areas. For example, it can be used for structural analysis of periodic descriptor systems \cite{Coll04,Verga07}. Also, we will encounter this equation in computing periodic deflating subspaces associated with a specified set of eigenvalues \cite{Granat07}. So, some scholars considered the numerical methods for computing the solution of the PGCS equation, see e.g., \cite{Chen12,Hajarian14} and references therein. It was also shown in \cite{Granat07} that if
\begin{equation*}\label{1.2}
\lambda(\{(A_k,C_k)\}_1^p)\cap\lambda(\{(B_k,D_k)\}_1^p)=\varnothing,
\end{equation*}
then the PGCS equation \eqref{1.1} has a unique solution. Here, $\lambda(\{(G_{k},H_{k})\}^{p}_{1})$ denotes the eigenvalue set of the periodic regular matrix pairs $\{(G_{k},H_{k})\}^{p}_{1}$ . This condition is equivalent to the fact that the coefficient matrix of the matrix-vector form of \eqref{1.1} is nonsingular. The matrix-vector form is
\begin{equation}\label{1.3}
Wz=g,
\end{equation}
where
$$W=\begin{bmatrix}
     I\otimes A_1 & -(B_1^T\otimes I) & 0 & 0 & 0 & \cdots & 0 & 0 \\
     0 & -(D_1^T\otimes I) & I\otimes C_1 & 0 & 0 & \cdots & 0 & 0 \\
     0 & 0 & I\otimes A_2 & -(B_2^T\otimes I) & 0 & \cdots & 0 & 0 \\
     0 & 0 & 0 & -(D_2^T\otimes I) &  I\otimes C_2 & \cdots & 0 & 0 \\
     \vdots & \vdots & \vdots & \vdots & \vdots & \ddots & \vdots & \vdots \\
     0 & 0 & 0 & 0 & 0 & \cdots & I\otimes A_p & -(B_p^T\otimes I) \\
     I\otimes C_p & 0 & 0 & 0 & 0 & \cdots & 0 & -(D_p^T\otimes I) \\
   \end{bmatrix},
$$
and
$$z={\rm vec}\left(\begin{bmatrix}
      X_1, Y_1,\cdots,X_p,Y_p
    \end{bmatrix}\right),\;
g={\rm vec}\left(\begin{bmatrix}
    E_1, F_1,\cdots,E_p,F_p
  \end{bmatrix}\right).
$$
In the above expressions, $ X\otimes Y$ denotes the Kronecker product \cite{Grah81}, the operator 'vec' stacks the columns of a matrix one underneath the other  \cite{Grah81}, $I$ is the identity matrix of appropriate order, and $K^T$ stands for the transpose of the matrix $K$.

For the similar motivations in \cite{Diao13,Diao12,Higham93,Kags94,Linw07}, we investigate the perturbation analysis for the PGCS equation in this paper. After introducing the notation and preliminaries in Section \ref{sec:N&P}, we present the normwise backward error for the PGCS equation in Section \ref{sec:B&E}. In Section \ref{sec:P&B}, the normwise and componentwise perturbation bounds for the PGCS equation are derived. A normwise condition number and the effective condition number are also given in this section. In Section \ref{sec:M&C}, we provide the mixed and componentwise condition numbers for the PGCS equation. An algorithm based on the SCE method is proposed to estimate the mixed and componentwise condition numbers in Section \ref{sec:N&E}. To estimate the normwise and effective condition numbers, we consider an alternative method, that is, the probabilistic spectral norm estimator by Hochstenbach \cite{Hochs13}, which provides a reliable estimation of the spectral norm. A corresponding algorithm is devised in Section \ref{sec:N&E}. In addition, the numerical examples are also given in this section to illustrate the differences between the normwise, effective, mixed and componentwise condition numbers, and the efficiency of the statistical condition estimations, respectively. Finally, we present the conclusion of the whole paper.

\section{Notation and preliminaries}\label{sec:N&P}
For the matrix $A=(a_{ij})\in \mathbb{R}^{m\times n}$, $A^\dag$, $\|A\|_2$, $\|A\|_{\infty}$, and $\|A\|_F$ stand for its Moore-Penrose inverse, spectral norm, max row norm, and Frobenius norm, respectively, $|A|$ is the matrix with elements $|a_{ij}|$, and $\left\| A \right\|_{\max }$ is defined by $\left\| A \right\|_{\max }=\left\| {{\rm vec}(A)} \right\|_\infty $. For the vectors $a=[a_1,\cdots,a_p]^T\in\mathbb{R}^{p}$ and $ b=[b_1,\cdots, b_p]^T \in\mathbb{R}^{p}$, we define the entry-wise division between $a$ and $b$ by ${a}/{b}=[c_{1},\cdots, c_{p}]^T$ with
$$c_i  = \left\{ \begin{array}{l}
 \frac{a_i}{b_i} , \textrm{ if }b_i  \ne 0, \\
 a_i,\textrm{ if } b_i  = 0. \\
 \end{array} \right. $$ Following \cite{R_Xie13}, the componentwise distance between $a$ and $b$ is defined by
\begin{equation*}
d(a,b) = {\left\| {\frac{{a - b}}{b}} \right\|_\infty } = \begin{array}{*{20}{c}}
{\mathop {\max }\limits_{i = 1, \cdots ,p} }
\end{array}\left\{ {\frac{{\left| {{a_i} - {b_i}} \right|}}{{\left| {{b_i}} \right|}}} \right\}= \left\{ \begin{array}{l}
 \frac{\left| a_{i_{0}} - b_{i_{0}}\right|}{\left| b_{i_{0}}\right|} , \textrm{ if }b_{i_{0}}  \ne 0, \\
 \left| a_{i_{0}}\right|, \quad \textrm{ if } b_{i_{0}}  = 0. \\
 \end{array} \right.
\end{equation*}
Note that when $b_{i_{0}}\neq 0$, $d(a,b)$ gives the relative distance from $a$ to $b$ with respect to $b$, while the absolute distance for $b_{i_{0}}= 0$.

In order to define the mixed and componentwise condition numbers, we also need to define the set $B^0(a,\epsilon)=\{x=[x_1,\cdots,x_p]^T\in \mathbb{R}^{p}\mid \left|x_i-a_i\right|\leq \epsilon |a_i|, i=1,\cdots,p\}$ with $a=[a_1,\cdots,a_p]^T\in \mathbb{R}^{p}$ and $\epsilon>0$, and denote the domain of definition of a function $F:{\mathbb{R}^p} \to {\mathbb{R}^q}$ by ${\rm Dom}(F)$. Thus, the definitions of the mixed and componentwise condition numbers can be given as follows.
\begin{definition}{\rm \cite{R_Xie13}}\label{Def M&Ccond}
Let $F:{\mathbb{R}^p} \to {\mathbb{R}^q}$ be a continuous map defined on an open set ${\rm Dom}(F) \subset {\mathbb{R}^p}$ such that $0 \notin {\mathop{\rm Dom}\nolimits}(F)$. Let $a \in {\rm Dom}(F)$, such that $F(a) \ne 0$.
\begin{enumerate}
  \item The mixed condition number of $F$ at $a$ is defined by
  \begin{equation*}
    m(F,a) = \mathop {\lim }\limits_{\epsilon  \to 0} \mathop {\sup }\limits_{\mathop {x \ne a}\limits_ {x \in {B^0}(a,\epsilon )}} \frac{{{{\left\| {F(x) - F(a)} \right\|}_\infty }}}{{{{\left\| {F(a)} \right\|}_\infty }}}\frac{1}{{d(x,a)}}.
  \end{equation*}

  \item The componentwise condition number of $F$ at $a$ is defined by
  \begin{equation*}
    c(F,a) = \mathop {\lim }\limits_{\epsilon  \to 0} \mathop {\sup }\limits_{\mathop {x \ne a}\limits_ {x \in {B^0}(a,\epsilon )}} \frac{{d(F(x),F(a))}}{{d(x,a)}}.
  \end{equation*}
\end{enumerate}
\end{definition}

The Fr\'{e}chet derivative is essential in deriving the explicit expressions of condition numbers. Its definition is presented below.
\begin{definition}\label{Def Frech diff}
 Let $\mathcal{X}$ and $\mathcal{Y}$ be two Banach spaces, and a map $F: U\in \mathcal{X} \rightarrow \mathcal{Y}$ with $U$ being an open set. Then $F$ is said to be Fr\'{e}chet differentiable at $a\in U$, if there exists a bounded linear operator $DF_a: \mathcal{X}\rightarrow \mathcal{Y}$ such that $$\lim_{h\rightarrow 0}\frac{\|F(a+h)-F(a)-DF_a(h)\|}{\|h\|}=0.$$
\end{definition}

When the map $F$ in Definition \ref{Def M&Ccond} is Fr\'{e}chet differentiable,  the following lemma given in \cite{R_Xie13} reduces the computation burden of mixed and componentwise condition numbers.
\begin{lem}\label{Lem expres m&c}
Under the assumptions of Definition \ref{Def M&Ccond}, when $F$ is Fr$\acute{e}$chet differentiable at $a$, we have
\begin{eqnarray}
  m(F,a) &=&\frac{\left\|\mid DF(a)\mid \mid a\mid \right\|_{\infty}}{\left\|F(a)\right\|_{\infty}} , \label{mix_cond}\\
   c(F,a)&=&\left\|\frac{\mid DF(a)\mid \mid a\mid}{\mid F(a)\mid}\right\|_{\infty},\label{comp_cond}
\end{eqnarray}
where $DF(a)$ is the Fr$\acute{e}$chet derivative of $F$ at $a$.
\end{lem}

To estimate the mixed and componentwise condition numbers, we need the SCE method which is ever mentioned in Section 1. In the following, we present a brief introduction on this method.

For a twice continuously differentiable function $f:\mathbb{R}^{p}\rightarrow \mathbb{R}$, by Taylor's theorem, we get
\begin{equation}\label{Taylor approx}
f(x+\delta z)=f(x)+\delta \nabla f(x)^Tz+O(\delta^2),
\end{equation}
where $\delta$ is a small positive number, $\nabla f(x)=\left[\frac{\partial f(x)}{\partial x_1}, \frac{\partial f(x)}{\partial x_2},\cdots,\frac{\partial f(x)}{\partial x_p}\right]^T$ is the derivative of $f$ at $x$, and $z\in \mathbb{R}^{p}$ satisfies $\|z\|_2=1$. From \eqref{Taylor approx}, the following inequality can be derived easily
$$|f(x+\delta z)-f(x)|\approx \delta |\nabla f(x)^Tz|\leqslant \delta\|\nabla f(x)\|_2,$$
which shows that the local sensitivity can be measured by a magnification factor $\delta$ and the absolute condition number $\|\nabla f(x)\|_2$. Based on the firm theoretical analysis given in \cite{Kenney94}, we have that if we choose a random vector $z$ from $\mathcal{U}(S_{p-1})$, the uniform distribution over unit sphere $S_{p-1}$ in $R^{p}$, then the following equality holds
\begin{equation}\label{eqn_Expect.Operat.}
\textbf{E}(\mid \nabla f(x)^Tz \mid)=\omega_p\|\nabla f(x)\|_2,
\end{equation}
where $\textbf{E}(\cdot)$ is the expectation operator, and $\omega_p$ is the Wallis factor with $\omega_1=1$, $\omega_2={2}/{\pi}$, and
$$\omega_p  =\left\{
  \begin{array}{ll}
    \frac{1\cdot 3 \cdot5 \cdots (p-2)}{2\cdot 4\cdot 6\cdots (p-1)}, & \hbox{for $p$ odd,} \\
    \frac{2}{\pi}\frac{2\cdot 4\cdot 6\cdots (p-2)}{3\cdot 5\cdot 7\cdots (p-1)}, & \hbox{for $p$ even,}
  \end{array}
\right.\textrm{ when } p>2.$$ Owing to the equality \eqref{eqn_Expect.Operat.} and the easy approximability of the Wallis factor ($\omega_p\approx \sqrt{\frac{2}{\pi(p-\frac{1}{2})}}$ preserves high accuracy),  $\eta={\mid \nabla f(x)^Tz\mid}/{\omega_p}$ can be used as a condition estimator, and satisfies the following probability relationship
\begin{equation*}
\textrm{Pr}(\frac{\|\nabla f(x)\|_2}{\gamma}\leq \eta \leq \gamma \|\nabla f(x)\|_2)\geqslant 1- \frac{2}{\pi \gamma}+O(\frac{1}{\gamma^2}),\quad{\rm with}\quad \gamma>1.
\end{equation*}
According to \cite{Kenney94}, the accuracy of condition estimator can be enhanced by multiple samples. If we choose two samples $\hat{z}_1$, $\hat{z}_2\in\mathcal{U}(S_{p-1})$, then the condition estimator given by
$$\eta(2)=\frac{\omega _2}{\omega_p}\sqrt{|\nabla f(x)^Tz_1|^2+ |\nabla f(x)^Tz_2|^2}$$
with $z_1$, $z_2$ being obtained from $\hat{z}_1$ and $\hat{z}_2$ by orthonormalization meets the following probability relationship
\begin{eqnarray*}
 {\rm Pr}(\frac{\|\nabla f(x)\|_2}{\gamma}\leq \eta(2) \leq \gamma \|\nabla f(x)\|_2)&\approx& 1- \frac{\pi}{4 \gamma^2}.
\end{eqnarray*}
In the similar manner, a general $k$-sample SCE estimator can be defined \cite{Kenney94}.

In addition, in the following sections, we will apply the following equality frequently
\begin{eqnarray}
{\rm vec}(AXC)& = &({C^T} \otimes A){\rm vec}(X),
\label{Kron1}
\end{eqnarray}
where $A,X$ and $C$ are matrices of appropriate orders such that the product $AXC$ is well-defined. The equality \eqref{Kron1} can be found in \cite{Grah81}.

\section{Normwise backward error}\label{sec:B&E}
Let $\mathcal{\widehat{Z}}=\left[\hat{X}_{1},\hat{Y}_{1},\cdots,\hat{X}_{p},\hat{Y}_{p} \right]$ denote an approximate solution to the PGCS equation \eqref{1.1}. The normwise backward error of $\mathcal{\widehat{Z}}$  is defined by
\begin{align}\label{def:backer}
&\eta(\mathcal{\widehat{Z}}) \equiv \min\{\epsilon: (A_k+\Delta A_k)\hat{X}_k-\hat{Y}_k(B_k+\Delta B_k)=E_k+\Delta E_k, \nonumber \\
&\quad\quad\quad\quad\quad\quad\quad (C_k+\Delta C_k)\hat{X}_{k+1}-\hat{Y}_k(D_k+\Delta D_k)=F_k+\Delta F_k,\; k=1,\cdots,p \},
\end{align}
where $\Delta A_{k},\Delta C_{k}\in\mathbb{R}^{m\times m}$, $\Delta B_{k},\Delta D_{k}\in\mathbb{R}^{n\times n}$, and $\Delta E_{k},\Delta F_{k}\in\mathbb{R}^{m\times n}$ satisfy
\begin{equation}\label{3.2}
\|\Delta A_k\|_F\leq \epsilon \alpha_k, \|\Delta B_k\|_F\leq \epsilon \beta_k, \|\Delta E_k\|_F\leq \epsilon \gamma_k, \|\Delta C_k\|_F\leq \epsilon \zeta_k, \|\Delta D_k\|_F\leq \epsilon \tau_k, \|\Delta F_k\|_F\leq \epsilon \delta_k.
\end{equation}
The tolerances $\alpha_{k},\beta_{k},\gamma_{k},\zeta_{k},\tau_{k}$  and $\delta_{k}$   provide some freedom in how we measure the perturbations. Usually,
\begin{equation}\label{3.200}
\alpha_{k}=\|A_{k}\|_{F},\ \beta_{k}=\|B_{k}\|_{F},\ \gamma_{k}=\|E_{k}\|_{F},\ \zeta_{k}=\|C_{k}\|_{F},\ \tau_{k}=\|D_{k}\|_{F},\ \delta_{k}=\|F_{k}\|_{F}.
\end{equation}
In this case, the normwise backward error is called the relative normwise backward error with respect to Frobenius norm.

The equation in \eqref{def:backer} can be rewritten as
\begin{equation}\label{3.3}
    \left\{
       \begin{array}{ll}
         \Delta A_k\hat{X}_k-\hat{Y}_k\Delta B_k-\Delta E_k=E_k-(A_k \hat{X}_k-\hat{Y}_k B_k)=R_{k1}, &  \\
          & \hbox{$k=1,\cdots,p,$} \\
         \Delta C_k\hat{X}_{k+1}-\hat{Y}_k\Delta D_k-\Delta F_k=F_k-(C_k \hat{X}_{k+1}-\hat{Y}_k D_k)=R_{k2}, & \hbox{}
       \end{array}
     \right.
\end{equation}
where $\mathcal{R}=\left[R_{11},R_{12},\cdots, R_{p1},R_{p2}\right]$ denotes the residual corresponding to the solution $\mathcal{\widehat{Z}}$.  Using the Kronecker product and \eqref{Kron1}, we can rewrite \eqref{3.3} as
\begin{equation*}
    \left\{
       \begin{array}{ll}
         (\hat{X}^T_k\otimes I){\rm vec}(\Delta A_k)-(I\otimes \hat{Y}_k){\rm vec}(\Delta B_k)-{\rm vec}(\Delta E_k)={\rm vec}(R_{k1}), &  \\
          & \hbox{$k=1,\cdots,p.$} \\
         (\hat{X}^T_{k+1}\otimes I){\rm vec}(\Delta C_k)-(I\otimes\hat{Y}_k){\rm vec}(\Delta D_k)-{\rm vec}(\Delta F_k)={\rm vec}(R_{k2}), & \hbox{}
       \end{array}
     \right.
\end{equation*}
That is,
\begin{equation}\label{3.4}
\widehat{H}u=r,
\end{equation}
where $\widehat{H}={\rm diag}(\widehat{H}_1,\cdots,\widehat{H}_p)$ with
\begin{align*}
&\widehat{H}_k = \left[ {\begin{array}{*{20}c}
                      \alpha_k(\hat{X}^T_k\otimes I) & -\beta_k(I\otimes\hat{Y}_k) & -\gamma_kI & 0 & 0 & 0 \\
                      0 & 0& 0 & \zeta_k(\hat{X}^T_{k+1}\otimes I) & -\tau_k(I\otimes\hat{Y}_k) & -\delta_kI  \\
                    \end{array}} \right],
\end{align*}
and
\begin{align*}
&u =\left[
           \frac{{\rm vec}(\Delta A_1)^T}{\alpha_1},\frac{{\rm vec}(\Delta B_1)^T}{\beta_1},\frac{{\rm vec}(\Delta E_1)^T}{\gamma_1},\frac{{\rm vec}(\Delta C_1)^T}{\zeta_1},\frac{{\rm vec}(\Delta D_1)^T}{\tau_1},\frac{{\rm vec}(\Delta F_1)^T}{\delta_1},\cdots,\right. \\
&\quad\quad\quad\left. \frac{{\rm vec}(\Delta A_p)^T}{\alpha_p},\frac{{\rm vec}(\Delta B_p)^T}{\beta_p},\frac{{\rm vec}(\Delta E_p)^T}{\gamma_p},\frac{{\rm vec}(\Delta C_p)^T}{\zeta_p},\frac{{\rm vec}(\Delta D_p)^T}{\tau_p},\frac{{\rm vec}(\Delta F_p)^T}{\delta_p}
         \right]^T,\\
&r={\rm vec}\left(\left[R_{11},R_{21},\cdots,R_{p1},R_{p2}\right]\right).
\end{align*}
It is easy to find that $\hat{H}$  is full row rank if $\gamma_{k}\neq0$ and $\delta_{k}\neq0$ for $k=1,\cdots,p.$  In this case, \eqref{3.4} has a minimum Euclidean norm solution
\begin{equation*}\label{3.5}
    u=\widehat{H}^{\dag}r.
\end{equation*}
From the definition of normwise backward error, we have
\begin{equation*}\label{3.6}
   \eta(\mathcal{\widehat{Z}})\leq \left\|\widehat{H}^{\dag}r\right\|_2.
\end{equation*}
On the other hand, considering \eqref{3.2},
\begin{equation*}\label{3.7}
   \|u\|_2^2=
\sum\limits_{i = 1}^p {\frac{{\left\| {\Delta A_i } \right\|_F^2 }}{{\alpha _i^2 }} + \frac{{\left\| {\Delta B_i } \right\|_F^2 }}{{\beta _i^2 }} + \frac{{\left\| {\Delta E_i } \right\|_F^2 }}{{\gamma _i^2 }} + \frac{{\left\| {\Delta C_i } \right\|_F^2 }}{{\zeta _i^2 }} + \frac{{\left\| {\Delta D_i } \right\|_F^2 }}{{\tau _i^2 }} + \frac{{\left\| {\Delta F_i } \right\|_F^2 }}{{\delta _i^2 }}}
\leq6p\epsilon^2.
\end{equation*}
Therefore,
\begin{equation}\label{3.8}
    \frac{\left\|\widehat{H}^{\dag}r\right\|_2}{\sqrt{6p}}\leq \eta(\mathcal{\widehat{Z}})\leq \left\|\widehat{H}^{\dag}r\right\|_2.
\end{equation}
Thus, we obtain both the upper and lower bounds of the normwise backward error $\eta(\mathcal{\widehat{Z}})$ for the PGCS equation.
\begin{rmk}
If the period $p=1$, the bounds in \eqref{3.8} reduce to the corresponding ones for the GCS equation. The reduced lower bound is a little different from the one in \cite{Kags94} since the definitions of normwise backward error here and in \cite{Kags94} are a little different. Further, if $C_1=0,D_1=0$, and $F_1=0$, we have the results for the classic Sylvester equation \cite{Higham93}. Note that $\sqrt{6}$ should be replaced by $\sqrt{3}$ in this case.
\end{rmk}

\section{Perturbation bounds}\label{sec:P&B}
Assume that the matrices $A_{k},B_{k},E_k,C_{k},D_{k},F_{k}, X_k$  and $Y_{k}$ in \eqref{1.1} are perturbed as
\begin{eqnarray*}
  &&A_k\rightarrow A_k+\Delta A_k,  B_k\rightarrow B_k+\Delta B_k,   E_k\rightarrow E_k+\Delta E_k,  X_k\rightarrow X_k+\Delta X_k, \\
   &&   C_k\rightarrow C_k+\Delta C_k, D_k\rightarrow D_k+\Delta D_k, F_k\rightarrow F_k+\Delta F_k, Y_k\rightarrow Y_k+\Delta Y_k,
\end{eqnarray*}
where $\Delta A_{k},\Delta C_{k}\in\mathbb{R}^{m\times m}$, $\Delta B_{k},\Delta D_{k}\in\mathbb{R}^{n\times n}$, $\Delta E_{k},\Delta F_{k},\Delta X_{k},\Delta Y_{k}\in\mathbb{R}^{m\times n}$, and $\Delta X_{p+1}=\Delta X_1$. Then the perturbed PGCS equation \eqref{1.1} is
\begin{equation}\label{4.1}
    \left\{
       \begin{array}{ll}
         (A_k+\Delta A_k)(X_k+\Delta X_k)-(Y_k+\Delta Y_k)(B_k+\Delta B_k)=E_k+\Delta E_k, & \hbox{} \\
          & \hbox{$k=1,\cdots,p.$} \\
         (C_k+\Delta C_k)(X_{k+1}+\Delta X_{k+1})-(Y_k+\Delta Y_k)(D_k+\Delta D_k)=F_k+\Delta F_k, & \hbox{}
       \end{array}
     \right.
\end{equation}
In the following, we regard $\Delta X_{k},\Delta Y_{k}(k=1,\cdots,p)$  as the unknown matrices of the matrix equation \eqref{4.1}, and obtain the condition under which the equation \eqref{4.1} has the unique solution, and then the desired perturbation bounds.

Considering \eqref{1.1}, the equation \eqref{4.1} can be simplified as
\begin{equation*}
    \left\{
       \begin{array}{ll}
         A_k\Delta X_k-\Delta Y_kB_k=\Delta E_k-(\Delta A_k X_k-Y_k\Delta B_k)-(\Delta A_k\Delta X_k-\Delta Y_k\Delta B_k), & \hbox{} \\
          & \hbox{$k=1,\cdots,p.$} \\
         C_k\Delta X_{k+1}-\Delta Y_kD_k=\Delta F_k-(\Delta C_k X_{k+1}-Y_k\Delta D_k)-(\Delta C_k\Delta X_{k+1}-\Delta Y_k\Delta D_k), & \hbox{}
       \end{array}
     \right.
\end{equation*}
which, using the Kronecker product and \eqref{Kron1}, can be rewritten as
\begin{equation}\label{4.2}
    W\begin{bmatrix}
       {\rm vec}(\Delta X_1) \\
       {\rm vec}(\Delta Y_1) \\
       \vdots \\
       {\rm vec}(\Delta X_p) \\
       {\rm vec}(\Delta Y_p) \\
     \end{bmatrix}=
     \begin{bmatrix}
       {\rm vec}(\Delta E_1) \\
       {\rm vec}(\Delta F_1) \\
       \vdots \\
       {\rm vec}(\Delta E_p) \\
       {\rm vec}(\Delta F_p) \\
     \end{bmatrix}-
     \Delta W\begin{bmatrix}
      {\rm vec}(X_1) \\
       {\rm vec}(Y_1) \\
      \vdots \\
       {\rm vec}(X_p) \\
       {\rm vec}(Y_p) \\
     \end{bmatrix}-
     \Delta W\begin{bmatrix}
       {\rm vec}(\Delta X_1) \\
       {\rm vec}(\Delta Y_1) \\
      \vdots \\
       {\rm vec}(\Delta X_p) \\
       {\rm vec}(\Delta Y_p) \\
     \end{bmatrix},
\end{equation}
where $\Delta W$ is the same as $W$ in \eqref{1.3} with  $A_{k},B_{k},C_{k}$, and $D_{k}$ being replaced by $\Delta A_{k},\Delta B_{k}, \Delta C_{k}$, and $\Delta D_{k},$
respectively. Let
\begin{equation*}
\Delta z={\rm vec}\left(\begin{bmatrix}
       \Delta X_1,\Delta Y_1,\cdots,\Delta X_p,\Delta Y_p
     \end{bmatrix}\right),\ \Delta g=
     {\rm vec}\left(\begin{bmatrix}
       \Delta E_1,\Delta F_1,\cdots,\Delta E_p,\Delta F_p
     \end{bmatrix}\right).
\end{equation*}
Then we simplify \eqref{4.2} as
\begin{equation}\label{4.3}
    W\Delta z=\Delta g-\Delta Wz-\Delta W\Delta z.
\end{equation}
Combining the first two terms in the right side of \eqref{4.3}, we can rewrite \eqref{4.3} as
\begin{equation}\label{4.4}
    W\Delta z=-{H}_1u-\Delta W\Delta z,
\end{equation}
where ${H}_1$ is the same as $\hat{H}$ in \eqref{3.4} except that $\hat{X}_{k}$ and $\hat{Y}_{k}$ in \eqref{3.4} are replaced by $X_{k}$ and $Y_{k}$, respectively. Thus,
\begin{align*}
    \Delta z=-W^{-1}{H}_1u-W^{-1}\Delta W\Delta z.
\end{align*}
Define the operator equation of $\Delta z $ as follows
\begin{equation}\label{4.5}
    \Phi(\Delta z)=\Delta z=-W^{-1}{H}_1u-W^{-1}\Delta W\Delta z.
\end{equation}
In the following, we use the Banach fixed point theorem (see,
e.g., \cite[Appendix D]{Kons03}) to derive the bound for $\Delta z $.

Let\begin{equation}\label{4.6}
   \left\|W^{-1}\Delta W\right\|_2<1,
\end{equation}
and denote the set $\Omega$ as
\begin{equation*}
    \Omega=\left\{s\in \mathbb{R}^{2mnp}:\|s\|_2\leq \frac{\left\|W^{-1}{H}_1u\right\|_2}{1-\left\|W^{-1}\Delta W\right\|_2}\right\},
\end{equation*}
which is closed and convex. Then, for any $s_{1},s_{2}\in\Omega$, we have
\begin{eqnarray*}
  \left\|\Phi(s_1)\right\|_2 &\leq& \left\|W^{-1}{H}_1u\right\|_2+\left\|W^{-1}\Delta W\right\|_2\|s_1\|_2 \\
  &\leq& \left\|W^{-1}{H}_1u\right\|_2+ \left\|W^{-1}\Delta W\right\|_2\frac{\left\|W^{-1}{H}_1u\right\|_2}{1-\left\|W^{-1}\Delta W\right\|_2}=\frac{\left\|W^{-1}{H}_1u\right\|_2}{1-\left\|W^{-1}\Delta W\right\|_2},
\end{eqnarray*}
and
\begin{equation*}
\left\|\Phi(s_1)-\Phi(s_2)\right\|_2\leq\left\|W^{-1}\Delta W s_1-W^{-1}\Delta W s_2\right\|_2\leq\left\|W^{-1}\Delta W\right\|_2\|s_1-s_2\|_2.
\end{equation*}
Therefore, $\Phi(\cdot)$  maps the set $\Omega$  into itself and is contractive (see,
e.g., \cite[Appendix D]{Kons03}). According to the Banach fixed point theorem, we have that there is a unique solution $\Delta z $  to the equation \eqref{4.5} in the set $\Omega$  when \eqref{4.6} holds.
As a result,
\begin{align}
&\|\Delta z\|_2 = \left\|[\Delta X,\Delta Y_1,\cdots, \Delta X_p,\Delta Y_p]\right\|_F\leq \frac{\left\|W^{-1}{H}_1u\right\|_2}{1-\left\|W^{-1}\Delta W\right\|_2}\nonumber\\
&\quad\quad\quad\leq \frac{\left\|W^{-1}{H}_1\right\|_2\left(
\sum\limits_{i = 1}^p {\frac{{\left\| {\Delta A_i } \right\|_F^2 }}{{\alpha _i^2 }} + \frac{{\left\| {\Delta B_i } \right\|_F^2 }}{{\beta _i^2 }} + \frac{{\left\| {\Delta E_i } \right\|_F^2 }}{{\gamma _i^2 }} + \frac{{\left\| {\Delta C_i } \right\|_F^2 }}{{\zeta _i^2 }} + \frac{{\left\| {\Delta D_i } \right\|_F^2 }}{{\tau _i^2 }} + \frac{{\left\| {\Delta F_i } \right\|_F^2 }}{{\delta _i^2 }}}\right)^{1/2}}{1-\left\|W^{-1}\Delta W\right\|_2}.\label{4.7}
\end{align}
What's more, if set
\begin{eqnarray*}\label{4.8}
&&\epsilon=\max\left\{\frac{\|\Delta A_1\|_F}{\alpha_1},\frac{\|\Delta B_1\|_F}{\beta_1},\frac{\|\Delta E_1\|_F}{\gamma_1},\frac{\|\Delta C_1\|_F}{\zeta_1},\frac{\|\Delta D_1\|_F}{\tau_1},\frac{\|\Delta F_1\|_F}{\delta_1},\cdots, \right. \nonumber\\
&&\quad\quad\quad\quad\quad\left. \frac{\|\Delta A_p\|_F}{\alpha_p},\frac{\|\Delta B_p\|_F}{\beta_p},\frac{\|\Delta E_p\|_F}{\gamma_p},\frac{\|\Delta C_p\|_F}{\zeta_p},\frac{\|\Delta D_p\|_F}{\tau_p}, \frac{\|\Delta F_p\|_F}{\delta_{p}}\right\},
\end{eqnarray*}
then we have
\begin{eqnarray}\label{4.9}
    \left\|[\Delta X,\Delta Y_1,\cdots, \Delta X_p,\Delta Y_p]\right\|_F&\leq&\frac{\sqrt{6p}\left\|W^{-1}{H}_1\right\|_2\epsilon}{1-\left\|W^{-1}\Delta W\right\|_2}.
\end{eqnarray}
In summary, we have the following theorem.
\begin{thm}\label{thm4.1}
Assume that the unperturbed and perturbed PGCS equations are given in \eqref{1.1} and \eqref{4.1}, respectively. If the perturbations in \eqref{4.1} satisfy \eqref{4.6}, then the perturbed PGCS equation \eqref{4.1} has a unique solution, and the normwise perturbation bounds \eqref{4.7} and \eqref{4.9} hold.
\end{thm}

\begin{rmk}
From \eqref{4.4} or \eqref{4.9}, by omitting the high-order terms, we can get the following first-order perturbation bound
\begin{eqnarray}
  \left\|[\Delta X,\Delta Y_1,\cdots, \Delta X_p,\Delta Y_p]\right\|_F\lesssim \sqrt{6 p}\left\|W^{-1}{H}_1\right\|_2 \epsilon.
\end{eqnarray}
The above bound is attainable to first-order in $\epsilon$. So,
\begin{eqnarray}\label{4.100}
k_{N1}=\frac{\left\|W^{-1}{H}_1\right\|_2}{\left\|[ X_1, Y_1,\cdots, X_p, Y_p]\right\|_F}
\end{eqnarray}
can be regarded as the normwise condition number for the PGCS equation \eqref{1.1}. It is a generalization of the ones for the GCS equation and the classic Sylvester equation given in \cite{Higham93,Kags94}.
\end{rmk}

\begin{rmk}
Using the equation \eqref{4.3}, along the same line for deriving \eqref{4.7}, we have the following bound under the condition \eqref{4.6},
\begin{eqnarray*}
  \frac{\left\|[\Delta X_1,\Delta Y_1,\cdots,\Delta X_p,\Delta Y_p]\right\|_F}{\left\|[ X_1, Y_1,\cdots, X_p, Y_p]\right\|_F} &\leq& \frac{\left[\frac{\|W^{-1}\|_2\|\Delta g\|_2}{\left\|[ X_1, Y_1,\cdots, X_p, Y_p]\right\|_F}+\left\|W^{-1}\Delta W\right\|_2 \right]}{1-\left\|W^{-1}\Delta W\right\|_2}  \\
   &\leq&  \frac{\left[\frac{\left\|[\Delta E_1,\Delta F_1,\cdots,\Delta E_p,\Delta F_p]\right\|_F}{\left\|[ E_1, F_1,\cdots, E_p, F_p]\right\|_F}k_{E} +\frac{\|\Delta W\|_2}{\|W\|_2}k(W)\right]}{1-\left\|W^{-1}\Delta W\right\|_2},
\end{eqnarray*}
where $k(W)=\|W\|_2\|W^{-1}\|_2$ and
\begin{eqnarray*}\label{4.101}
k_{E}=\frac{\|W^{-1}\|_2\left\|[ E_1, F_1,\cdots, E_p, F_p]\right\|_F}{\left\|[ X_1, Y_1,\cdots, X_p, Y_p]\right\|_F}.
\end{eqnarray*}
As done in \cite{Diao13}, we can call $k_{E}$ the effective condition number for the PGCS equation \eqref{1.1}. It can be much tighter than $k_{N1}$ if there are only perturbations on the right-hand side of the equation \eqref{1.1}. The main reason is that $k_{E}$ only contains the information of $[ E_1, F_1,\cdots, E_p, F_p]$, while $k_{N1}$ contains the information of all the coefficient matrices.
\end{rmk}

Now we consider the componentwise perturbation bounds for the PGCS equation using the operator equation \eqref{4.5} and the generalized Banach fixed point theorem (see,
e.g., \cite[Appendix D]{Kons03}).

Let
\begin{equation}\label{4.10}
   {\rm radius}\left(\left|W^{-1}\Delta W\right|\right)<1,
\end{equation}
and define the set $\Xi$ as
$$\Xi=\left\{s\in \mathbb{R}^{mnp}: |s|\leq (I-|W^{-1}\Delta W)|)^{-1}|W^{-1}H_1u|\right\}.$$
It is easy to check that the set $\Xi$ is closed and convex, and for any $s_1,s_2\in \Xi$,
\begin{eqnarray*}
  |\Phi(s_1)| &\leq& \left|W^{-1}H_1u\right|+\left|W^{-1}\Delta W\right||s_1| \\
   &\leq&  \left|W^{-1}H_1u\right|+\left|W^{-1}\Delta W\right|(I-\left|W^{-1}\Delta W)\right|)^{-1}\left|W^{-1}H_1u\right| \\
   &=&(I-|W^{-1}\Delta W)|)^{-1}|W^{-1}H_1u|
\end{eqnarray*}
and
$$|\Phi(s_1)-\Phi(s_2)|\leq \left|W^{-1}\Delta Ws_1-W^{-1}\Delta Ws_2\right|\leq \left|W^{-1}\Delta W\right||s_1-s_2|.$$
Therefore, $\Phi(v,\cdot)$ maps the set $\Xi$ into itself and is generalized contractive (see,
e.g., \cite[Appendix D]{Kons03}). According to the generalized Banach fixed point theorem, we have that there is a unique solution $\Delta z$ to the equation \eqref{4.5} in the set $\Xi$  when \eqref{4.10} is satisfied. As a result,
\begin{equation}\label{4.11}
    |\Delta z|={\rm vec}\left(\begin{bmatrix}
                 \left|\Delta X_1,\Delta Y_1,\cdots,\Delta X_p,\Delta Y_p\right|
               \end{bmatrix}\right)\leq(I-\left|W^{-1}\Delta W)\right|)^{-1}\left|W^{-1}H_1u\right|.
\end{equation}
The above discussions imply the following theorem.

\begin{thm}\label{thm4.2}
Assume that the unperturbed and perturbed PGCS equations are given in \eqref{1.1} and \eqref{4.1}, respectively. If the perturbations in \eqref{4.1} fulfill \eqref{4.10}, then the perturbed PGCS equation \eqref{4.1} has a unique solution, and the componentwise perturbation bound \eqref{4.11} holds.
\end{thm}

\begin{rmk}\label{rmk4.3}
Form \eqref{4.11}, we have the first-order componentwise perturbation bound
\begin{equation}\label{4.12}
    {\rm vec}\left(\begin{bmatrix}
                 \left|\Delta X_1,\Delta Y_1,\cdots,\Delta X_p,\Delta Y_p\right|
               \end{bmatrix}\right)\lesssim\left|W^{-1}H_1u\right|.
\end{equation}
\end{rmk}
\begin{rmk}
When the period $p=1$, the perturbation bounds obtained in this section reduce to the corresponding ones for the GCS equation, where the first-order normwise one is equivalent to the one in \cite{Kags94} in essence.
\end{rmk}

\section{Mixed and componentwise condition numbers}\label{sec:M&C}
In this section, using Lemma 2.1, we investigate the mixed and componentwise condition numbers for the PGCS equation

We first rewrite \eqref{4.4}, omitting the high-order terms, as follows,
$$W\Delta z \approx -{H}_2v,$$
where $H_2$ and $v$ are the same as $H_1$ and $u$ in \eqref{4.4}, respectively, except that all the tolerances $\alpha_k,\beta_k,\gamma_k,\zeta_k,\tau_k$ and $\delta_k$ are replaced by 1. Thus,
\begin{equation}\label{5.1}
    \Delta z\approx -W^{-1}{H}_2v.
\end{equation}
Define the map $\Psi$ as
\begin{equation*}
    \Psi:t\rightarrow z,
\end{equation*}
where $t=$$[{\rm vec}(A_1)^T$, ${\rm vec}(B_1)^T$, ${\rm vec}(E_1)^T$, ${\rm vec}(C_1)^T$, ${\rm vec}(D_1)^T$, ${\rm vec}(F_1)^T$,$\cdots$, ${\rm vec}(A_p)^T$, ${\rm vec}(B_p)^T$, ${\rm vec}(E_p)^T$, ${\rm vec}(C_p)^T$, ${\rm vec}(D_p)^T$, ${\rm vec}(E_p)^T]^T$, and $z$ is defined as in \eqref{1.3}. Then from Definition \ref{Def Frech diff} and \eqref{5.1}, it follows that the Fr\'{e}chet derivative of $\Psi$  at $t$ is:
\begin{equation}\label{5.2}
    D\Psi(t)=-W^{-1}{H}_2.
\end{equation}
Thus, combining Lemma \ref{Lem expres m&c} with \eqref{5.2}, we have the following theorem which gives the expressions of the mixed and componentwise condition numbers of the PGCS equation \eqref{1.1}.

\begin{thm}\label{thm5.1}
With the above notation, the mixed and componentwise condition numbers of the PGCS equation \eqref{1.1} are given by
\begin{eqnarray}\label{5.3}
&&m(\Psi,t)=\frac{\left\|\left|W^{-1}{{H}_2}\right||t|\right\|_{\infty}}{\left\|[X_1, Y_1,\cdots, X_p, Y_p]\right\|_{\max}}=\frac{\left\|\omega\right\|_{\infty}}{\left\|[X_1, Y_1,\cdots, X_p, Y_p]\right\|_{\max}},\\
\label{5.4}
&&c(\Psi,t)=\left\|\frac{\left|W^{-1}{{H}_2}\right||t|}{{\rm vec}([X_1, Y_1,\cdots, X_p, Y_p])}\right\|_{\infty}=\left\|\frac{\omega}{{\rm vec}([X_1, Y_1,\cdots, X_p, Y_p])}\right\|_{\infty},
\end{eqnarray}
where
\begin{eqnarray*}
  \omega &=&  \left|W^{-1}\begin{bmatrix}
                    X_1^T\otimes I \\
                    0 \\
                  \end{bmatrix}
   \right||{\rm vec}(A_1)|+\left|W^{-1}\begin{bmatrix}
                    I\otimes Y_1 \\
                    0 \\
                  \end{bmatrix}
   \right||{\rm vec}(B_1)|+\left|W^{-1}\begin{bmatrix}
                     I \\
                    0 \\
                  \end{bmatrix}
   \right||{\rm vec}(E_1)|  \\
   && +\left|W^{-1}\begin{bmatrix}
                    0 \\
                    X_2^T\otimes I \\
                    0 \\
                  \end{bmatrix}
   \right||{\rm vec}(C_1)|+\left|W^{-1}\begin{bmatrix}
                    0 \\
                    I\otimes Y_1 \\
                    0 \\
                  \end{bmatrix}
   \right||{\rm vec}(D_1)|+\left|W^{-1}\begin{bmatrix}
                     0 \\
                    I \\
                    0 \\
                  \end{bmatrix}
   \right||{\rm vec}(F_1)| \\
   & &+\cdots +\left|W^{-1}\begin{bmatrix}
                    0 \\
                    X_p^T\otimes I \\
                    0\\
                  \end{bmatrix}
   \right||{\rm vec}(A_p)|+\left|W^{-1}\begin{bmatrix}
                    0 \\
                    I\otimes Y_p \\
                    0\\
                  \end{bmatrix}
   \right||{\rm vec}(B_p)|+\left|W^{-1}\begin{bmatrix}
                     0 \\
                    I \\
                    0\\
                  \end{bmatrix}
   \right||{\rm vec}(E_p)|   \\
   &&  +\left|W^{-1}\begin{bmatrix}
                     0\\
                    X_1^T\otimes I \\
                  \end{bmatrix}
   \right||{\rm vec}(C_p)|+\left|W^{-1}\begin{bmatrix}
                    0 \\
                    I\otimes Y_p\\
                  \end{bmatrix}
   \right||{\rm vec}(D_p)| +\left|W^{-1}\begin{bmatrix}
                     0 \\
                    I\\
                  \end{bmatrix}
   \right||{\rm vec}(E_p)|.
\end{eqnarray*}
\end{thm}
\begin{pf}
In view of Lemma \ref{Lem expres m&c}, \eqref{5.2} and the definition of $\left\|\cdot\right\|_{\max}$, it is necessary only to show how to obtain the expression of $\omega$. This can be done easily by using the expressions of $H_2$ and $t$.  $\square$
\end{pf}

\begin{rmk}
Note that
\begin{equation*}\label{5.5}
    \|\omega\|_{\infty}=\left\|\left|W^{-1}{{H}_2}\right||t|\right\|_{\infty}\leq \left\|W^{-1}\right\|_{\infty}\left\|\left|{H}_2\right||t|\right\|_{\infty}=\left\|W^{-1}\right\|_{\infty} \left\|\begin{bmatrix}
             |A_1||X_1|+|Y_1||B_1|+|E_1| \\
             |C_1||X_2|+|Y_1||D_1|+|F_1| \\
             \vdots \\
             |A_p||X_p|+|Y_p||B_p|+|E_p| \\
             |C_p||X_1|+|Y_p||D_p|+|F_p| \\
           \end{bmatrix}
    \right\|_{\max}.
\end{equation*}
Here, the definition of $\left\|\cdot\right\|_{\max}$ is used. So, we have an upper bound for the mixed condition number
\begin{equation}\label{5.6}
    m(\Psi,t)\leq \frac{\left\|W^{-1}\right\|_{\infty}}{\left\|[X_1, Y_1,\cdots, X_p, Y_p]\right\|_{\max}}\left\|\begin{bmatrix}
             |A_1||X_1|+|Y_1||B_1|+|E_1| \\
             |C_1||X_2|+|Y_1||D_1|+|F_1| \\
             \vdots \\
             |A_p||X_p|+|Y_p||B_p|+|E_p| \\
             |C_p||X_1|+|Y_p||D_p|+|F_p| \\
           \end{bmatrix}
    \right\|_{\max}.
\end{equation}
From the definition of the entry-wise division of vectors given in Section 2, we have
\begin{equation*}
\left\|\frac{\left|W^{-1}{{H}_2}\right||t|}{{\rm vec}([X_1, Y_1,\cdots, X_p, Y_p])}\right\|_{\infty}= \left\|\left({\rm diag}\left({\rm vec}([X_1, Y_1,\cdots, X_p, Y_p])\right)\right)^{\ddag}\left|W^{-1}{{H}_2}\right||t|\right\|_{\infty}.
\end{equation*}
Here, for a vector $a=[a_1,\cdots,a_p]^T\in\mathbb{R}^{p}$, $\left({\rm diag}(a)\right)^{\ddag}$ denotes a diagonal matrix with the elements $a_i^{\ddag} (i=1,\cdots,p)$ of the following form
$$a_i^{\ddag}  = \left\{ \begin{array}{l}
 \frac{1}{a_i} , \textrm{ if }a_i  \ne 0, \\
 1,\ \textrm{ if } a_i  = 0. \\
 \end{array} \right. $$
Then
\begin{equation*}
    c(\Psi, t)\leq \left\|\left({\rm diag}\left({\rm vec}([X_1, Y_1,\cdots, X_p, Y_p])\right)\right)^{\ddag}W^{-1}\right\|_{\infty}\left\|\left|{H}_2\right||t|\right\|_{\infty}.
\end{equation*}
Thus, an upper bound for the componentwise condition number can be given by
\begin{equation}\label{5.7}
    c(\Psi, t)\leq \left\|\left({\rm diag}\left({\rm vec}([X_1, Y_1,\cdots, X_p, Y_p])\right)\right)^{\ddag}W^{-1}\right\|_{\infty}\left\|\begin{bmatrix}
             |A_1||X_1|+|Y_1||B_1|+|E_1| \\
             |C_1||X_2|+|Y_1||D_1|+|F_1| \\
             \vdots \\
             |A_p||X_p|+|Y_p||B_p|+|E_p| \\
             |C_p||X_1|+|Y_p||D_p|+|F_p| \\
           \end{bmatrix}
    \right\|_{\max}.
\end{equation}
\end{rmk}

\begin{rmk}
Using \eqref{5.2} and the definition of the normwise condition number given in \cite{Rice66}, we can obtain an alternative normwise condition number for the PGCS equation \eqref{1.1}:
\begin{eqnarray*}\label{5.100}
k_{N2}=\frac{\left\|W^{-1}{H}_2\right\|_2\left(
\sum\limits_{i = 1}^p\left\|A_i\right\|_F^2+ \left\|B_i\right\|_F^2+ \left\|E_i\right\|_F^2+ \left\|C_i\right\|_F^2+ \left\|D_i\right\|_F^2+ \left\|F_i\right\|_F^2\right)^{1/2}}{\left\|[ X_1, Y_1,\cdots, X_p, Y_p]\right\|_F},
\end{eqnarray*}
which is a little larger than $k_{N1}$ in \eqref{4.100} if the tolerances $\alpha_k,\beta_k,\gamma_k,\zeta_k,\tau_k$ and $\delta_k$ in \eqref{4.100} are chosen as in \eqref{3.200}. In addition, if the period $p=1$, the above condition number reduces to the corresponding one for the GCS equation \cite{Linw07}.
\end{rmk}

\section{Numerical experiments}\label{sec:N&E}
In this part, our attention mainly focuses on the comparison and estimation of the condition numbers derived in the above sections.

We first provide an example to compare the normwise, effective, mixed and componentwise condition numbers. The example is taken from \cite{Chen12} with some modifications.
\begin{example}\label{example1}{\rm For the PGCS equation \eqref{1.1}, let the period $p=3$, and the coefficient matrices be
\begin{align*}
& A_1=\begin{bmatrix}
       1 & 0 & 0.1 \\
       0 & 1 & 10 \\
       0 & 0 & 1 \\
     \end{bmatrix}, A_2=\begin{bmatrix}
          1 & 0.3 & 8 \\
          0 & 1 & 10 \\
          0 & 0 & 1 \\
        \end{bmatrix}, A_3=\begin{bmatrix}
          0.1 & 0.03 & 9 \\
          0 & 0.1 & 0.9 \\
          0 & 0 & 0.1 \\
        \end{bmatrix},
       B_1=\begin{bmatrix}
          1 & 12 \\
          0 & 2\\
        \end{bmatrix},\\
& B_2=\begin{bmatrix}
          2 & 1 \\
          0 & 1\\
        \end{bmatrix}
, B_3=\begin{bmatrix}
          1 & 21 \\
          0 & 10^{-t}\\
        \end{bmatrix}, E_1=\begin{bmatrix}
         1 & 1 \\
         0 & 1 \\
         0 & 10 \\
       \end{bmatrix}
, E_2=\begin{bmatrix}
         0 & 1 \\
         2 & 1 \\
         5 & 8 \\
       \end{bmatrix}
, E_3=\begin{bmatrix}
         2 & 0 \\
         3 & 1 \\
         0 & 2 \\
       \end{bmatrix},\\
& C_1=\begin{bmatrix}
       0.1 & 10 & 1.5 \\
       1 & 10 & 0.1 \\
       2 & 0.3 & 0.1 \\
     \end{bmatrix}
, C_2=\begin{bmatrix}
       1.1 & 3 & 8 \\
       0.2 & 5 & 0.1 \\
       1 & 0.01 & 0.01 \\
     \end{bmatrix}
, C_3=\begin{bmatrix}
       1 & 0.5 & 0.9 \\
       1 & 0.1 & 0.9 \\
       1 & 2 & 0.15 \\
     \end{bmatrix}
, D_1=\begin{bmatrix}
          1 & 0 \\
          1 & 2 \\
        \end{bmatrix},\\
& D_2=\begin{bmatrix}
          2 & 9 \\
          2 & 1 \\
        \end{bmatrix}
, D_3=\begin{bmatrix}
          1 & 1 \\
          3 & 10^{-\tau} \\
        \end{bmatrix}
, F_1=\begin{bmatrix}
          1 & 0 \\
          0.1 & 1 \\
          2 & 0\\
        \end{bmatrix}
, F_2=\begin{bmatrix}
          0 & 1 \\
          2 & 1 \\
          5 & 8 \\
        \end{bmatrix}
, F_3=\begin{bmatrix}
          2 & 0 \\
          3 & 1 \\
          2 & 5 \\
        \end{bmatrix},\\
& \mathrm{with}\; \tau,\; t\; \in \{1,3,5\}.
\end{align*}
Upon some computations, the numerical results are exhibited in Table I.
\begin{table}[htbp]\label{Table1}
\centering
 \begin{tabular}{lrrrrrr}
  \toprule
   \multicolumn{7}{c}{\centering  Table I: Comparison of condition numbers } \\
\cmidrule(r){1-7}
& $\tau=1$ & $\tau=1$&$\tau=1$ & $\tau=3$& $\tau=3$  & $\tau=5$ \ \\
& $t=1$ &$t=3$&$t=5$ & $t=3$& $t=5$  &$t=5$ \ \\
\midrule
$k_{N1}$  & 564.1934 &  1.4085e+003 & 1.3455e+003 &  1.4065e+003 & 1.3438e+003 & 1.3438e+003\\
$k_{N2}$ &2.3429e+004 &  5.8489e+004 & 5.5874e+004 &  5.8407e+004  &  5.5803e+004 & 5.5803e+004 \\
$k_{E}$  & 263.9046 &  182.1415 &  181.5541 &    182.1423 & 181.5566 & 181.5567 \\
$ m(\Psi,t)$ & 52.9059  &  18.1312 & 16.1057  &   18.1240 &   16.1058 &   16.1058 \\
$c(\Psi,t)$  & 1.3318e+003 & 260.1651 & 269.9788 & 120.0864  &   119.9581 & 119.9582\\
\bottomrule
\end{tabular}
\end{table}
}\end{example}
%\begin{rmk}
From Table I, one can easily find that the effective, mixed and componentwise condition numbers behave well in most cases, while the normwise condition numbers $k_{N1}$ and $k_{N2}$ may highly overestimate the condition of the PGCS equation. Here, it should be pointed out that $c(\Psi,t)$ may be very large if there are very small elements in the solution. This may be the reason why $c(\Psi,t)$ is so large for $\tau=1$ and $t=1$. In this case, some distinction should be made to cope with this extremal case. We suggest the projection method proposed by Arioli et al. \cite{Arioli07}, and  Cao and Petzold \cite{Cao03}, but we will not go that far in this paper.
%\end{rmk}

In the following, we will devise two algorithms based on the probabilistic spectral norm estimator and the SCE method to estimate the normwise, effective, mixed and componentwise condition numbers. The former will be called the PCE method for short.

\begin{algorithm}[htbp]{\small
\caption{PCE for the normwise and effective condition numbers}\label{AlgorithmPCE}
\begin{enumerate}
  \item Generate a starting vector $v_0$ from $\mathcal{U}(S_{q-1})$ with $q=2p(m^2+n^2+mn)$.
  \item Compute the guaranteed lower bound $\alpha$ and the probabilistic upper bound $\beta$ of $\left\|W^{-1}{H}_1\right\|_2$ ($\|W^{-1}\|_2$) by probabilistic spectral norm estimator ($\left\|W^{-1}{H}_1\right\|_2\leqslant \beta$ ($\|W^{-1}\|_2\leqslant \beta$) will hold with a given probability $1-\epsilon$, where $\epsilon$ is a user-chosen parameter).
  \item Compute the normwise and effective condition number by
  \begin{equation*}
  k_{pceN1}=\frac{\alpha+\beta}{2\left\|[ X_1, Y_1,\cdots, X_p, Y_p]\right\|_F},\quad \left(k_{pceE}= \frac{\alpha+\beta}{2}\frac{\left\|( E_1, F_1,\cdots, E_p, F_p)\right\|_F}{\left\|( X_1, Y_1,\cdots, X_p, Y_p)\right\|_F}\right).
  \end{equation*}
\end{enumerate}}
\end{algorithm}

\begin{algorithm}[htbp]{\small
\caption{SCE for the mixed and componentwise condition numbers}\label{AlgorithmSCE}
\begin{enumerate}
  \item Generate the random matrices $(R_{11},L_{11},M_{11},S_{11},N_{11},Q_{11},\cdots,R_{p1}, L_{p1}, M_{p1}, S_{p1},N_{p1},Q_{p1})$, $\cdots,$ $(R_{1s}, L_{1s},M_{1s},S_{1s},N_{1s},Q_{1s},\cdots, R_{ps}, L_{ps}, M_{ps},S_{ps},N_{ps},Q_{ps})$, where $R_{kj}, S_{kj}\in\mathbb{R}^{m\times m}$, $L_{kj},N_{kj}\in\mathbb{R}^{n\times n}$ and $M_{kj},Q_{kj}\in\mathbb{R}^{m\times n}$ with $k=1,\cdots,p$, $j=1,\cdots,s$, and all entries being in the standard normal distribution $\mathcal{N}(0,1)$. Orthonormalize the matrix
      $$\begin{bmatrix}
        {\rm vec}(R_{11}) & \cdots & {\rm vec}(R_{1s}) \\
        {\rm vec}(L_{11}) & \cdots & {\rm vec}(L_{1s}) \\
       \vdots & \cdots & \vdots \\
        {\rm vec}(Q_{p1}) & \cdots & {\rm vec}(Q_{ps}) \\
      \end{bmatrix}$$ to get an orthonormal matrix $[p_1,\cdots,p_s]$. Then, convert $p_j$ into the matrix form $$(R_{1j},L_{1j},M_{1j},S_{1j},N_{1j},Q_{1j},\cdots,R_{pj}, L_{pj}, M_{pj}, S_{pj},N_{pj},Q_{pj}).$$
  \item Set $q=2p(m^2+n^2+mn)$, get the approximates of $\omega_q$ and $\omega_s$, and let
  \begin{align*}
  &(R_{1j},L_{1j},M_{1j},S_{1j},N_{1j},Q_{1j},\cdots,R_{pj}, L_{pj}, M_{pj}, S_{pj},N_{pj},Q_{pj})\\
  &=(R_{1j},L_{1j},M_{1j},S_{1j},N_{1j},Q_{1j},\cdots,R_{pj}, L_{pj}, M_{pj}, S_{pj},N_{pj},Q_{pj})\\
  &\quad\circ ( A_1,B_1,  E_1, C_1,D_1,F_1\cdots,A_p,B_p,E_p,C_p, D_p,F_p).
  \end{align*}
     Here, the symbol $\circ$ denotes the Hadamard product.
  \item For $j=1,\cdots,s$, solve the following PGCS equation
  \begin{equation*}
    \left\{
       \begin{array}{ll}
         A_k X_{kj}-Y_{kj}B_k=M_{kj}-(R_{kj} X_k-Y_kL_{kj}), & \hbox{} \\
          & \hbox{$k=1,\cdots,p,$} \\
         C_k X_{(k+1)j}-Y_{kj}D_k=Q_{kj}-(S_{kj} X_{k+1}-Y_kN_{kj}), & \hbox{}
       \end{array}
     \right.
\end{equation*}
and compute the absolute condition vector
$$\kappa_{abs}=\frac{\omega_s}{\omega_q}\sqrt{\sum_{j=1}^{s}\mid u_j\mid^2},$$
where $u_j={\rm vec}\left(\begin{bmatrix}
      X_{1j}, Y_{1j},\cdots,X_{pj},Y_{pj}
    \end{bmatrix}\right)$. Here, the operations of taking square root and power are componentwise.
  \item Compute the estimations of the mixed and componentwise condition numbers by
  \begin{equation*}
  m_{sce}(\Psi,t)=\frac{\|\kappa_{abs}\|_{\infty}}{\left\|[X_1, Y_1,\cdots, X_p, Y_p]\right\|_{\max}}, \quad c_{sce}(\Psi,t)=\left\|\frac{\kappa_{abs}}{{\rm vec}([X_1, Y_1,\cdots, X_p, Y_p])}\right\|_{\infty}.
  \end{equation*}
\end{enumerate}
\emph{Note: For the sake of convenience, we write $( A_1,B_1,  E_1, C_1,D_1,F_1\cdots,A_p,B_p,E_p,C_p, D_p,F_p)$ as a matrix though the matrices in the parenthesis do not have same orders.}}
\end{algorithm}

The main part of Algorithm \ref{AlgorithmPCE} is to estimate $\left\|W^{-1}{H}_1\right\|_2$ ($\|W^{-1}\|_2$) by probabilistic spectral norm estimator. A detailed analysis of the estimator was given in \cite{Hochs13} by Hochstenbach. The author showed that $\left\|W^{-1}{H}_1\right\|_2$ ($\|W^{-1}\|_2$) can be contained in a small interval $[\alpha , \beta]$ with high probability. Here $\beta/ \alpha\leqslant 1+\delta$, where $\delta$ is another user-chosen parameter. In our computation, we take $\epsilon=0.001$ and $\delta= 0.01$ . Thus, $\left\|W^{-1}{H}_1\right\|_2\leqslant \beta$ ($\|W^{-1}\|_2\leqslant \beta$) holds with a probability at least $99.9\%$ and $\beta/ \alpha\leqslant 1.01$. Hence, we take $(\alpha+\beta)/2$ as the estimation of $\left\|W^{-1}{H}_1\right\|_2$ ($\|W^{-1}\|_2$).

For Algorithm \ref{AlgorithmSCE}, we would like to choose $s=3$ in numerical experiments. This means that $m_{sce}(\Psi,t)$ and $c_{sce}(\Psi,t)$ fall into the intervals $[0.2\times m(\Psi,t), 5\times m(\Psi,t)]$ and $[0.2\times c(\Psi,t), 5\times c(\Psi,t)]$ with the probability $1- \frac{32}{3\pi^2 \gamma^3}\approx 0.9913$, respectively, if $\gamma=5$.

Now we present a specific example to investigate the efficiency of these two algorithms in estimating the condition numbers.
\begin{example}\label{example2}
For the PGCS equation \eqref{1.1}, let $p=3$, $m=5$, and $n=4$, and generate the coefficient matrices as follows: $A_k, C_k\in \mathrm{\texttt{randn}}(m)$, $B_k, D_k\in \mathrm{\texttt{randn}}(n)$, and $E_k, F_k\in \mathrm{\texttt{randn}}(m,n)$. Here, the Matlab functions are used. Since the orders of the coefficient matrices are not so large, we get the solution by solving the linear equation \eqref{1.3}. The computed solution $\hat{z}$ satisfies the inequality $\||W^{-1}||r|\|_{\infty}/\|\hat{z}\|_{\infty}\leq 10^{-8}$ \cite[p.131]{R_Higham02} and is treated as the exact solution. We test $1000$ PGCS equations, and define the ratios of the estimated condition numbers and the exact ones as follows
\begin{align*}
 r_{N1}=\frac{k_{pceN1}}{k_{N1}},\quad r_E=\frac{k_{pceE}}{k_{E}}, \quad r_m=\frac{m_{sce}(\Psi,t)}{m(\Psi,t)},\quad r_c=\frac{c_{sce}(\Psi,t)}{c(\Psi,t)}.
\end{align*}
Upon computation, we have the numerical results of these ratios and their means and variances: $\mathbf{E}(r_{N1}) =1.0003$, $\mathbf{V}(r_{N1}) =5.7960e-007$, $\mathbf{E}(r_{E}) =1.0004$, $\mathbf{V}(r_{E}) =8.0694e-007$, $\mathbf{E}(r_{m}) =1.8313$, $\mathbf{V}(r_{m}) =2.4788$, $\mathbf{E}(r_{c}) =2.4269$, $\mathbf{V}(r_{c}) =7.1857$. The numerical results are plotted in Figure \ref{Figure1}.
\begin{figure}[htbp]
  \centering
  \includegraphics[width=1\textwidth,height=0.6\textwidth]{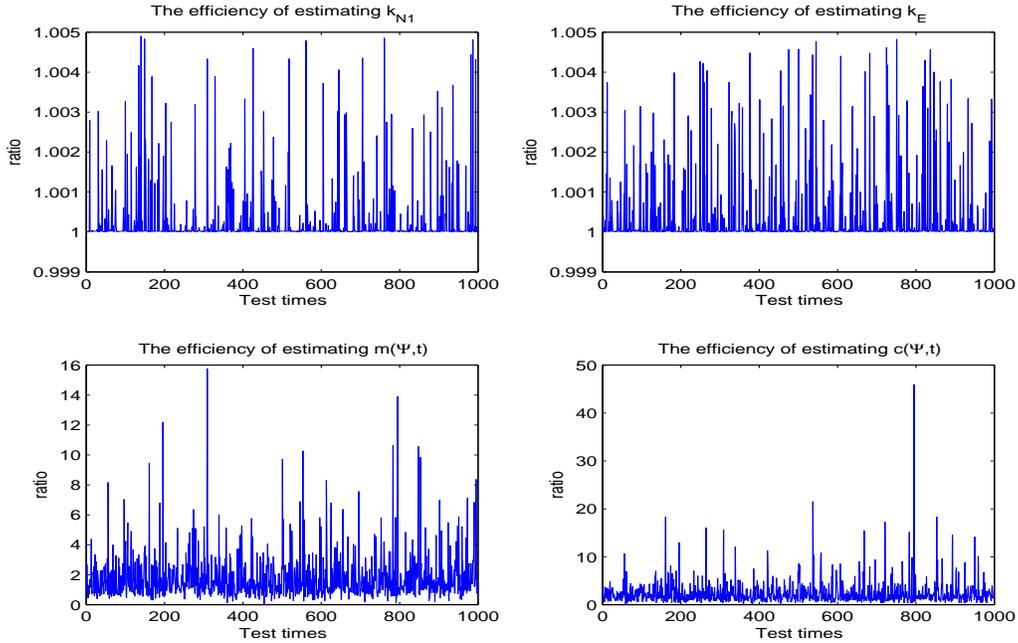}\\
  \caption{Efficiency of condition estimators}\label{Figure1}
\end{figure}
\end{example}
From Figure \ref{Figure1} and the results on means and variances, we can find that both the PCE method and the SCE method can give reliable estimations of the normwise, effective, mixed and componentwise condition numbers, respectively.
\begin{rmk}
In Example \ref{example2}, we get the solution to the PGCS equation \eqref{1.1} by solving the linear system \eqref{1.3}. The cost will be very expensive when the orders of the coefficient matrices in the PGCS equation \eqref{1.1} are large. In this case, other iterative methods need to be consulted; see \cite{Chen12,Hajarian14} and references therein.
\end{rmk}

\section{Conclusion}\label{sec:CONCL}
In this paper, we investigated the perturbation analysis of the PGCS equation. The normwise backward error for this equation is first given. Then, by Banach fixed point theorem, we derive its rigorous normwise and componentwise perturbation bounds, from which the first-order perturbation bounds, and the normwise and effective condition numbers  are obtained. Moreover, the explicit expressions of the mixed and componentwise condition numbers and their upper bounds for the PGCS equation are also given. A simple example is provided to illustrate the differences among these condition numbers. To estimate these condition numbers, the probabilistic spectral norm estimator and the SCE method are introduced and two algorithms are devised. From the numerical experiments, we find that both the PCE method and the SCE method perform efficiently in estimating the normwise, effective, mixed and componentwise condition numbers, respectively.\\

\ack{The work is supported by the National Natural Science Foundation of China under grant number 11201507 and the China Scholarship Council. The authors would like to thank Prof. Michiel E. Hochstenbach for providing Matlab program of probabilistic spectral norm estimator.} %This work was done while the first author was a visiting scholar at
%the Department of Mathematics and Statistics, Auburn University, from August 2014 to August
%2015.}


\begin{thebibliography}{99}
\bibitem{Arioli07}M. Arioli, M. Baboulin, S. Gratton, A partial condition number for linear least squares problems, SIAM J. Matrix Anal. Appl. 29(2) (2007) 413--433.
\bibitem{Cao03} Y. Cao, L. Petzold,  A subspace error estimate for linear systems, SIAM J. Matrix Anal. Appl. 24(3) (2003) 787--801.
\bibitem{Chen12}X. Chen, Solving the (generalized) periodic sylvester equation with the matrix sign function, Math. Numer. Sin. 34(2) (2012) 153--162 (in Chinese).
\bibitem{Coll04}C. Coll, M. Fullana, E. Sanchez, Reachability and observability indices of a discrete-time periodic descriptor system, Appl. Math. Comput. 153 (2004) 485--496.
%\bibitem{R_Cucker07}F. Cucker, H. Diao, Y. Wei, On mixed and componentwise condition numbers for Moore-Penrose inverse and linear least squares problems, Math. Comp. 76 (2007) 947--963.
\bibitem{Datta04} B. Datta, Numerical Methods for Linear Control Systems: Design and Analysis, Elsevier, London, 2003.
\bibitem{Demmel87}J. Demmel, B. K{\aa}gstr\"{o}m, Computing stable eigendecompositions of matrix pencils, Linear Algebra Appl. 88/89 (1987) 139--186.
\bibitem{Diao13}H. Diao, X. Shi, Y. Wei,  Effective condition numbers and small sample statistical condition estimation for the generalized Sylvester equation, Sci China Math 56 (2013) 967--982.
\bibitem{Diao12}H. Diao, H. Xiang, Y. Wei, Mixed, componentwise condition numbers and small sample statistical condition estimation of Sylvester equations, Numer. Linear Algebra Appl. 19 (2012) 639--654.
\bibitem{Ding05}F. Ding, T. Chen, Iterative least-squares solutions of coupled Sylvester matrix equations, Systems Control Lett. 54 (2005) 95--107.
\bibitem{Gohb93}I. Gohberg, I. Koltracht, Mixed, componentwise, and structured condition numbers, SIAM J. Matrix Anal. Appl. 14 (1993) 688--704.
\bibitem{Grah81}A. Graham, Kronecker Products and Matrix Calculus: with Applications, John Wiley, New York, 1981.
%\bibitem{Granat06a}R. Granat, B. K{\aa}gstr\"{o}m, Direct eigenvalue reordering in a product of matrices in periodic Schur form, SIAM J. Matrix Anal. Appl. 28 (2006) 285--300
\bibitem{Granat07}R. Granat, B. K{\aa}gstr\"{o}m, D. Kressner, Computing periodic deflating subspaces associated with a specified set of eigenvalues, BIT 47 (2007) 763--791.
\bibitem{Gudmun97}T. Gudmundsson, C. Kenney, A. Laub, Small-sample statistical estimates for the sensitivity of eigenvalue problems, SIAM J. Matrix Anal. Appl. 18 (1997) 868--886.
\bibitem{Hajarian14}M. Hajarian, Developing CGNE algorithm for the periodic discrete-time generalized coupled Sylvester matrix equations, Comp. Appl. Math. (2014) 1--17. Doi:10.1007/s40314-014-0138-7.
%\bibitem{Hench94}J. Hench, A. Laub, Numerical solution of the discrete-time periodic Riccati equation, IEEE Trans. Auto. Con. 39 (1994) 1197--1210.
\bibitem{R_Higham02}N. Higham, Accuracy and Stability of Numerical Algorithms, second ed., SIAM, Philadelphia, 2002.
%\bibitem{R_Higham08}N. Higham, Functions of Matrices: Theory and Computation, SIAM, Philadelphia, 2008.
\bibitem{Higham93}N. Higham, Perturbation theory and backward error for $AX-XB=C$, BIT  33 (1993) 124--136.
\bibitem{Hochs13}M. Hochstenbach, Probabilistic upper bounds for the matrix two-norm, J. Sci. Comput. 57 (2013) 464--476.
\bibitem{Jons02a}I. Jonsson, B. K{\aa}gstr\"{o}m, Recursive blocked algorithms for solving triangular systems-Part I: One-sided and coupled Sylvester-type matrix equations, ACM Trans. Math. Software 28 (2002) 392--415.
\bibitem{Jons02b}I. Jonsson, B. K{\aa}gstr\"{o}m, Recursive blocked algorithms for solving triangular systems-Part II: Two-sided and generalized Sylvester and Lyapunov matrix equations, ACM Trans. Math. Software  28 (2002) 416--435.
\bibitem{Kags94}B. K{\aa}gstr\"{o}m, A perturbation analysis of the generalized Sylvester equation $(AR -LB, DR -LE) = (C, F)$, SIAM J. Matrix Anal. Appl. 15 (1994) 1045--1060.
%\bibitem{Kaga96}B. K{\aa}gstr\"{o}m, P. Poromaa, LAPACK-style algorithms and software for solving the generalized Sylvester equation and estimating the separation between regular matrix pairs, ACM Trans. Math. Software 22 (1996) 78--103.
%\bibitem{Kags89}B. K{\aa}gstr\"{o}m, L. Westin, Generalized Schur methods with estimators for solving the generalized Sylvester equation, IEEE Trans. Automat. Control 34 (1989) 745--751.
\bibitem{Kenney94}C. Kenney, A. Laub, Small-sample statistical condition estimates for general matrix functions, SIAM J. Sci. Comput. 15 (1994) 36--61.
\bibitem{Kenney98a}C. Kenney, A. Laub, M. Reese, Statistical condition estimation for linear systems, SIAM J. Sci. Comput. 19 (1998) 566--583.
\bibitem{Keney98b}C. Kenney, A. Laub, M. Reese, Statistical condition estimation for linear least squares, SIAM J. Matrix Anal. Appl. 19 (1998) 906--923.
\bibitem{Kons03}M. Konstantinov, D. Gu, V. Mehrmann, P. Petkov, Perturbation Theory for Matrix Equations, Elsevier, Amsterdam, 2003.
\bibitem{laub08}A. Laub, J. Xia, Applications of statistical condition estimation to the solution of linear systems, Numer. Linear Algebra Appl. 15 (2008) 489--513.
\bibitem{Linw07}Y. Lin, Y. Wei, Condition numbers of the generalized Sylvester equation, IEEE Trans. Automat. Control 52 (2007) 2380--2385.
\bibitem{Rice66}J. Rice, A theory of condition, SIAM J. Numer. Anal. 3 (1966) 287--310.
\bibitem{Verga07}A .Varga, On computing minimal realizations of periodic descriptor systems. In: Proceedings of IFAC workshop on periodic control systems, St. Petersburg, Russia, 2007.
%\bibitem{Varga97}A. Varga, Periodic Lyapunov equations: some applications and new algotithms, Inter. J. Con 67 (1997) 69--87.
\bibitem{R_Xie13}Z. Xie, W. Li, X. Jin, On condition numbers for the canonical generalized polar decomposition of real matrices, Electron. J. Linear Algebra 26 (2013) 842--857.
\end{thebibliography}
\end{document}